\magnification 1200

\newcount \version
\newcount \shortversion
\newcount \longversion
\shortversion 1\relax
\longversion 2\relax
\version \longversion\relax


\font\smallrm=cmr8

\font\Bbb=msbm10
\def\BBB#1{\hbox{\Bbb#1}}

\font\Frak=eufm10
\def\frak#1{{\hbox{\Frak#1}}}

\def\gen{1.1}
\def\gentwo{1.2}
\def\eone{1.3}
\def\conv{1.4}
\def\nopt{1.5}
\def\obv{1.6}
\def\four{2.1}
\def\RXX{2.2}
\def\noptw{2.3}
\def\noptx{2.4}
\def\triop{2.5}
\def\trioptic{2.6}
\def\spi{2.7}
\def\Tone{2.8}
\def\RXY{2.9}
\def\RXZ{2.10}
\def\Ttwo{2.11}
\def\Tthree{2.12}
\def\Tfour{2.13}
\def\Exa{2.14}

\def\fdec{0.1}
\def\prob{1.1}
\def\CH{1.2}
\def\prec{1.3}
\def\betterprec{1.4}
\def\cost{1.5}
\def\wop{1.6}
\def\fourw{2.1}
\def\tanc{2.2}
\def\fourc{2.3}
\def\XYX{2.4}
\def\fsys{2.5}
\def\tang{2.6}
\def\tae{2.7}
\def\tax{2.8}
\def\eXYX{2.9}
\def\eYXY{2.10}
\def\ftau{2.11}
\def\fmu{2.12}
\def\XY{2.13}
\def\eqa{2.14}
\def\eqb{2.15}
\def\eqc{2.16}
\def\onef{2.17}

\def\g{{\frak g}}
\def\R{{\BBB R}}
\def\N{{\BBB N}}
\def\C{{\BBB C}}
\def\F{{\BBB F}}
\def\B{{\cal B}}
\def\bS{\overline{S}}
\def\tmin{t_{\hbox{\smallrm inf}}}
\def\H{{\BBB H}}
\def\ep{\epsilon}
\def\ad{\hbox{\rm ad}}
\def\min{\hbox{\rm min}}
\def\tx{t_x}
\def\ty{t_y}

\centerline
{\bf Time-optimal decompositions in $SU(2)$.}

\centerline
{ {\bf Yuly Billig} \kern-3pt
\footnote{${}^1$}{School of Mathematics and Statistics, Carleton University,
1125 Colonel By Drive, Ottawa, K1S 5B6, Canada. E-mail: billig@math.carleton.ca}}
\footnote{}{2010 Mathematics Subject Classification: 81P68, 22E70, 57R27.}

\

\

{\bf Abstract.}
{A connected Lie group $G$ is generated by its two 1-parametric subgroups $\exp(t X)$, 
$\exp(tY)$ if and only if the Lie algebra of $G$ is generated by $\{ X, Y \}$.
We consider decompositions of elements of $G$ into a product of such exponentials with times $t > 0$ and study the problem of minimizing the total time
of the decompositions for a fixed element of $G$. We solve this problem for 
the group $SU_2$ and describe the structure of the time-optimal decompositions.}

\

\

{\bf 0. Introduction.}

In the axiomatics of quantum mechanics a state of a quantum system is a unit vector in a complex 
hermitian vector space and its evolution is given by unitary transformations. In quantum 
computing the space of states is finite-dimensional and the quantum computation is an element 
of $SU_N$, which is a compact connected real Lie group. In order to carry out a quantum 
computation in physical reality, we need to be able to transform an initial state $\psi$ 
of a quantum processor into a state $g \psi$, where $g$ is the element of $SU_N$ 
representing the quantum computation. In the classical (non-quantum) case, a computation can be viewed 
as a boolean function in several variables. A processor can not implement an arbitrary 
function directly, but instead decomposes the computation into elementary steps, which 
is reflected in the process of programming. An efficient algorithm for a classical computation
is a program that minimizes the number of elementary steps, which results in a shorter run time.

Likewise, an implementation of a quantum computation $g$ is a factorization of $g$ into certain
elementary factors, called the quantum gates (in analogy with the classical logical gates used
for the decomposition of an arbitrary boolean function of several variables in a disjunctive normal
form, for example). Some quantum gates are the direct analogues of the classical logical gates and
are discrete, while others depend on a continuous time parameter. 

Time evolution of a quantum system
is governed by a Hamiltonian $H$ and is given by the exponential $\exp(itH) \in SU_N$.
Here $X = iH$ belongs to the Lie algebra $su_N$. We are going to assume that all available
quantum gates are continuous, since the discrete gates are physically realized as continuous gates
applied for a specific finite time.

This leads to the following quantum control problem: given a set $S \subset su_N$ of quantum
controls (gates), decompose an element $g\in SU_N$ into a product
$$ g = \exp(t_1 C_1) \times \ldots \times \exp(t_n C_n), \eqno{(\fdec)}$$
where $C_i \in S$, $t_i \in \R$. This quantum control problem admits a solution for every $g$
precisely when the set $S$ generates the Lie algebra $su_N$ (see Theorem \gen \ below). This result is proved using topological methods and does not provide an effective procedure
for finding such decompositions.

Minimal sets of quantum gates for quantum computing were discussed in [2] and [6],
where it is shown that a set of two generic Hamiltonians is sufficient for controllability. A review of quantum control methods in physical chemistry is given in [7].

 The problem of finding explicit factorizations of type (\fdec) goes back to Euler [3], 
who studied it for the group $SO_3$ of rotations in $\R^3$. The Lie algebra $so_3$ consists 
of skew-symmetric $3\times 3$ matrices and the exponentials $\exp(t C_i)$ corresponding to
$$C_1 = \pmatrix{0 & 0 & 0 \cr 0 & 0 & -1 \cr 0 & 1 & 0 \cr}, \quad 
C_2 = \pmatrix{0 & 0 & -1 \cr 0 & 0 & 0 \cr 1 & 0 & 0 \cr}$$
are the rotations in angle $t$ around $X$- and $Y$-axis respectively. Euler proved that 
every $g \in SO_3$ can be factored as $g = \exp(t_1 C_1) \exp(t_2 C_2) \exp(t_3 C_1)$
for some $t_1, t_2, t_3 \in [0, 2\pi)$. The parameters $t_1, t_2, t_3$ are called Euler angles. 
This problem has applications in robotics when an object needs to be oriented in space in a 
prescribed way.

For the quantum applications the parameters $t_i$ in (\fdec) represent time, and thus must be positive. 
It is natural to pose an optimization problem in this context: among all factorizations (\fdec)
with positive times $t_i$, find the one that minimizes the total time $t_1 + \ldots + t_n$.
This corresponds to finding an implementation of a quantum algorithm with the shortest running time.

As we shall see in this paper, this optimization problem, as posed, might not have a solution -- there
may be a sequence of factorizations of type (\fdec) with the number of factors $n$ going to infinity,
while the total time going to infinum. For this reason we slightly modify the above optimization
problem and ask for the infinum of total times for all factorizations of a given $g\in SU_N$

In this paper we solve this problem for the group $SU_2$ with the set of controls $S$ consisting of
two elements, $S = \{ X, Y \}$. 

We show that the infinum time does not change if we replace $S = \{ X, Y \}$ with its convex closure
$\bS = \{ \tau X + (1-\tau) Y \, | \, 0 \leq \tau \leq 1 \}$. In the latter case there will be in fact an 
optimal decomposition (\fdec) with a finite number of factors. It turns out that to get an optimal
decomposition it is sufficient to add to $S$ at most one element $W \in \bS$, the one that 
is orthogonal to $X - Y$. In this paper we give explicit descriptions of the
optimal decompositions in $SU_2$.
Since there is a surjective homomorphism $SU_2 \rightarrow SO_3$, our results are also applicable to the group $SO_3$.

Our solution is based on the method of Lagrange multipliers adopted for the set-up of Lie groups.
Alternatively one could use the geometric theory developed in [4], which is based on the Pontryagin
maximum principle. The general methods, however, give only necessary conditions for optimality,
which in practice are not sufficient. To get the desired results we supplement these general methods
with some explicit calculations in $SU_2$ which allow us to get stronger and more explicit 
optimality conditions.

We hope that our methods and results will help to solve this problem in greater generality.

\

{\bf Acknowledgements:} I am thankful to Velimir Jurdjevic and Alexander Weekes for stimulating discussions. This work is supported with a grant from the Natural Sciences and Engineering Research Council of Canada.
 
\

\

{\bf 1. Time-optimal decompositions in Lie groups.}

Let $G$ be a compact connected real Lie group, and let $\g$ be its Lie algebra.
An element $X \in\g$ defines a 1-parametric subgroup
$\{ \exp (tX) | t \in \R \}$. In this paper we study an optimal control problem on $G$,
describing optimal decompositions of an arbitrary given element $g \in G$ into 
a product of exponentials $\exp(tX)$ with $X$ belonging to a fixed set $S$ of controls, 
$S \subset \g$ and positive times $t$.

The following well-known criterion describes when the group $G$ is controllable by a set
$S$:

{\bf Theorem \gen.  ([4], Theorem 6.1)} A real connected Lie group $G$ is generated by its subgroups 
$\{ \exp (t X) \}$, $X \in S$, if and only if $S$ generates the Lie algebra $\g$. 


In Theorem \gen \ the parameters in the subgroups $\{ \exp (t X) \}$ could be both positive 
and negative. In applications to controllability of quantum systems, the parameters $t$ 
represent time, and must be positive. If we restrict the question to controllability with 
positive time parameters, then the analogue of Theorem \gen \ still holds for compact Lie groups:

{\bf Theorem \gentwo. ([4], Theorem 6.3)} Let $G$ be a compact connected real group.
If the set $S$ generates the Lie algebra of $G$ then every element $g\in G$ admits
a factorization 
$$ g = \exp(t_1 C_1) \cdot \ldots \cdot \exp(t_n C_n)$$
for some $n \geq 0$ with $C_i \in S$ and $t_i > 0$.

These theorems are proved using topological methods, and do not provide an effective way 
of finding such decompositions. In this context it is natural to pose the problem of describing decompositions that are time-optimal: 

{\bf Problem.} For a given $g \in G$ determine
$$\inf \big\{ \, t_1 + \ldots + t_n \,  \big|  \, 
g = \exp(t_1 C_1) \cdot \ldots \cdot \exp( t_n C_n),
\; t_i \geq 0, C_i \in S  \, \big\}.  \eqno{(\prob)}$$


A compact Lie group $G$ is isomorphic to a Lie subgroup in a general linear group ([5], Corollary 4.22), so we assume that $G \subset GL_d (\F)$, where $\F = \R$ or $\C$.
The Lie algebra $\g$ is then a real subalgebra in the matrix Lie algebra $M_d (\F)$.

Consider an arbitrary $\R$-bilinear positive-definite scalar product $\left<\cdot , \cdot \right>$ on $M_d (\F)$.
We can make it left and right $G$-invariant by integration over $G$ using the Haar measure on $G$
([5], Section IV.2):
$$\left( A , B \right) =
 \int\limits_{{\kern 5pt}G \times G}{\kern -13pt}\int \left< gAh , gBh\right> dg dh .$$
Since the Haar measure on $G$ is left- and right-invariant, this averaging procedure will yield a $G$-bi-invariant scalar product $(\cdot \, ,\cdot)$ on $M_d (\F)$ which is still positive-definite. This scalar product defines a norm $| \cdot |$ on $M_d (\F)$.
We can rescale the above scalar product to achieve
$$|AB| \leq |A| \cdot |B| \quad \hbox{\rm for all \ } A, B \in M_d (\F).$$ 
The restriction of this scalar product to $\g$ gives a positive-definite bilinear form on $\g$, which is invariant under the conjugation  action of $G$ on its Lie algebra:
$$( g A g^{-1} , \, g B g^{-1}) = (A, B) \quad \hbox{\rm for  \ } A, B \in \g, \; g\in G .$$
Note that the Killing form on $\g$ yields a proportional scalar product when $\g$ is simple, but in case of a general compact Lie group $G$ the Killing form has a disadvantage that it is only (negative) semi-definite ([5], Corollary 4.26), and may be degenerate, whereas the form constructed above is strictly positive-definite.

We get the following estimate on the norm of a Lie bracket:
$$| [A, B] | = | AB - BA | \leq 2 \cdot |A| \cdot |B| .$$

The exponential map 
$$\exp: \; \g \rightarrow G, \quad \exp(A) = \sum\limits_{k = 0}^\infty {1\over k!} A^k$$
is surjective for compact $G$ ([5], Corollary 4.48), but is not injective. However locally, it maps bijectively a neighbourhood of $0\in \g$ to a neighbourhood of $1\in G$.
The logarithm, which is a local inverse of the exponential map and assumes values in $\g$, is defined for 
$g\in G$ satisfying $|1 - g| < 1$ using the standard formula
$$\ln (g) = - \sum\limits_{k=1}^\infty {1\over k} B^k , \quad \hbox{\rm where \ }
B = 1 - g.$$

The following Lemma is a straightforward norm estimation with power series of operators, and we omit its proof:

{\bf Lemma \eone.} (a) $| \exp (A) | \leq e^{|A|} $,

(b)  $|\exp(A+B) - \exp(A) | \leq  e^{|A|} \left(  e^{|B|} - 1 \right)$, 

(c)  $|\exp(A+B) - \exp(A) | \leq  e^{|A|} |B| + o(|B|)$, as $|B| \to 0$.

%
%
%

\

 Let us recall the Campbell-Hausdorff formula ([1], II.6.4):
$$\exp(A) \exp(B) = \exp \left( \sum\limits_{r,s \geq 0} H_{r,s} (A,B) \right), 
\eqno{(\CH)}$$
where
$$H_{r,s} (A,B) = {1\over r+s} \sum\limits_{m\geq 1} { (-1)^{m-1} \over m}
\sum\limits_{{r_1 + \ldots r_m = r \atop s_1 + \ldots s_m = s} \atop r_i + s_i \geq 1}
\left( \prod\limits_{i=1}^m {1 \over r_i! s_i!} \right)
[A^{r_1} B^{s_1} \ldots A^{r_m} B^{s_m}]$$
and $H_{0,0} (A,B) = 0$.
Here for a monomial $U = A^{r_1} B^{s_1} \ldots A^{r_m} B^{s_m}$, we denote
by $[U]$ the iterated commutator, defined inductively for $C \in \{ A, B \}$ as
$[C] = C$, $[C U] = [C, [U]]$. The iterated commutator 
$ [A^{r_1} B^{s_1} \ldots A^{r_m} B^{s_m}]$ is zero unless either
$s_m = 1$ or $s_m = 0$ and $r_m =1$.

The first few terms of this series are:
$$ \sum\limits_{r,s \geq 0} H_{r,s} (A,B)  = A + B + {1 \over 2} [A,B] 
+ {1\over 12} [A, [A,B]] - {1\over 12} [B, [A,B]] + \ldots .$$ 

Let us discuss convergence of the Campbell-Hausdorff series 
$\sum\limits_{r,s \geq 0} H_{r,s} (A,B) $. 
Consider a function in two variables
$$f(u,v) = - {1\over 2} \ln ( 2 - e^ {2u+2v} ) .$$
We can expand this function in a Taylor series
$$ f(u,v) =\sum\limits_{r,s \geq 0} \eta_{r,s} u^r v^s ,$$
which is absolutely convergent when $|u| + |v| \leq {\ln 2 \over 2}$.

By Lemma II.7.1 in [1], we can estimate the terms of the Campbell-Hausdorff series:
$$| H_{r,s} (A,B) | \leq \eta_{r,s} |A|^r |B|^s .$$
Thus the Campbell-Hausdorff series is convergent in $\g$ and the equality (\CH) holds
when $|A| + |B| \leq {\ln 2 \over 2}$. In addition, a tail of the Campbell-Hausdorff series
may be estimated with the corresponding tail of the series $f( |A|, |B|)$.

 We will also use the following well-known formula ([5], Proposition 1.93):
$$\exp (A) B \exp(-A) = \exp( \ad A ) B ,$$
where $\ad (A) B = [A, B]$.

{\bf Theorem \conv.}  Optimization problem (\prob) in $SU_2$ with $S = \{ X, Y \}$ is equivalent to the problem with the set of controls $\bS =  \{ \tau X + (1-\tau)Y \,| \tau \in [0,1] \}$.

{\bf Proof.} We need to show that for every $\tau \in [0,1]$ the infinum of time in (\prob)
with $S = \{ X, Y \}$ for $g = \exp ( \tau X + (1-\tau)Y)$ is less or equal to $1$. It is well-known that $g$ can be approximated with an arbitrary precision in time $1$
(see e.g., Theorem 3.7 in [4]).
 Indeed,
let $A = \tau X$, $B = (1-\tau)Y$. Then
$$ g = \exp (A+B) = 
\lim\limits_{N \to \infty} \left[ \exp \left( {A \over N} \right) \exp \left( {B \over N} \right) \right]^N . \eqno{(\prec)}$$
To see this, we apply the Campbell-Hausdorff formula,
$$ \exp \left( {A \over N} \right) \exp \left( {B \over N} \right) 
= \exp \left( {A \over N} + {B \over N} + {1 \over 2N^2} [A, B] + o({1/ N^2}) 
\right),$$
from which (\prec) follows: 
$$ \left[ \exp \left( {A \over N} \right) \exp \left( {B \over N} \right) \right]^N
= \exp( A + B + O(1/N)) .$$

We need to prove a stronger result, that $g$ itself can be decomposed in a product of exponentials with control set $S$ with a total time not exceeding $1+\epsilon$ for an arbitrarily small $\epsilon > 0$.


A direct computation with the Campbell-Hausdorff formula shows that we can improve the
approximation (\prec) in the following way: 
$$  h_N = \left[ \exp \left( {A \over 2N} \right) \exp \left( {B \over N} \right)  \exp \left( {A \over 2N} \right) \right]^N
= \exp\left( A + B + O(1/N^2) \right) . \eqno{(\betterprec)}$$
Let us show that we can attain $\exp \left( A+B \right)$ by a small variation of time parameters in (\betterprec), while keeping these parameters positive.
Let $\B(\epsilon) = \{ g\in G \, \big| \, |1-g| \leq \epsilon \}$.
Since the norm is $G$-invariant, we get that the ball in $G$ of radius $\epsilon$ with the center in $h\in G$ can be written as $h \B(\epsilon) = \B(\epsilon) h$.
We are going to show that by varying three of the parameters by $N^{- {3\over 2}}$, we can cover the set
$ h_N \B\left( c N^{- {3\over 2}} \right)$, which will contain $\exp(A+B)$ for large $N$ by Lemma \eone (c).
In order to show that a variation of certain three parameters covers a ball around $h_N$, we need to establish that the corresponding jacobian is non-zero. 
The computation of the jacobian takes place in the Lie algebra $su_2$. Since $\{ A, B \}$ is a generating set for this Lie algebra, the set $\{ A, B, [A,B] \}$ forms a basis of $su_2$. The tangent vector to $\exp \left( t \cdot \ad(A+B) \right) B$ at $t = 0$ is $[A,B]$, which implies that for some $M \in \N$, 
the Lie algebra $su_2$ is spanned by $\left\{ A, B , \exp \left( M^{-1} \ad (A+B) \right) B \right\}$.
Fix $M$, and assume that $N$ is a multiple of $M$, $N = MK$. We consider the following variation of (\betterprec):
$$\exp\left(  {A \over 2N} + \ep_1  A \right)
\exp\left( {B \over N} + \ep_2 B \right)
\left[ \exp\left( {A \over N} \right)\exp\left( {B \over N} \right) \right]^{K-1} \exp\left( {A \over N} \right)$$
$$ \times \exp\left( {B \over N} + \ep_3 B \right)
\left[ \exp\left( {A \over N} \right)\exp\left( {B \over N} \right) \right]^{N-K-1} \exp\left( {A \over 2N} \right) .$$
The differential of this variation is 
$$\exp\left( {A \over 2N} \right)
\left\{ \ep_1 A + \ep_2 B + \ep_3 \left[\exp\left( {B \over N} \right)\exp\left( {A \over N} \right) \right]^K B \left[\exp\left( {B \over N} \right)\exp\left( {A \over N} \right) \right]^{-K} \right\}$$
$$\times \left[ \exp\left( {B \over N} \right)\exp\left( {A \over N} \right) \right]^{N-1}
\exp\left( {B \over N} \right)\exp\left( {A \over 2N} \right)  .$$
Since 
$$\lim\limits_{N \to \infty}  \left[\exp\left( {B \over N} \right)\exp\left( {A \over N} \right) \right]^K B  \left[\exp\left( {B \over N} \right)\exp\left( {A \over N} \right) \right]^{-K}=  
\exp \left( {M^{-1}} \ad (A+B) \right) B ,$$
we conclude that for large $N$ the jacobian of this variation is non-zero. Thus by varying the corresponding time parameters with
$|\ep_1|, |\ep_2|, |\ep_3| \leq N^{- {3\over 2}}$, we can cover the set $ h_N  \B\left( c N^{- {3\over 2}} \right)$, for some constant $c>0$ that depends on $A$ and $B$.
Since $h_N$ approximates $g = \exp (A+B)$ with precision $O(N^{-2})$, this set contains $g$ for large $N$. This completes the proof of the Theorem.

Since a uniform rescaling of the controls $X^\prime = c X$, $Y^\prime = cY$, gives an equivalent optimization problem (with optimal time rescaled by a factor $c^{-1}$), we may assume without 
loss of generality that  $|X| = 1$ and $|Y| \geq |X|$. Let $\kappa = {1 / |Y|} \leq 1$.
We may pass to the normalized set of controls $\{ X, {Y / |Y|} \}$ by replacing
the total time $t_1 + \ldots + t_n$ for the decomposition
$\exp(t_1 C_1) \cdot \ldots \cdot \exp(t_n C_n)$ with the cost function
$$\sum_{i=1}^n \kappa_i t_i , \quad \hbox{\rm where} \;
\kappa_i = \left\{ \matrix{&1, &\hbox{\rm if \ } C_i  = X, \cr
&\kappa, &\hbox{\rm if \ } C_i  = Y/|Y|. \cr} \right. \eqno{(\cost)}$$

 From now on we assume that $S = \{ X, Y \}$ with $|X| = |Y| = 1$ and
consider the optimization problem with cost function (\cost) where 
the {\it cost factor} $\kappa \leq 1$.
This will allow us to consider the limiting case $\kappa = 0$, when there is no cost
associated with control $Y$.

Let us introduce some terminology. 

An {\it admissible word of length $n$} is an expression
$\exp(t_1 C_1) \cdot \ldots \cdot \exp(t_n C_n)$ with $n \geq 0$, $t_i \geq 0$ and 
$C_i \in S$.  A word of zero length is the identity element of $G$.

Every admissible word can be written in a  {\it reduced form}, where
$t_i > 0$ and $C_i \neq C_{i+1}$ for all $i$. 

A decomposition of $g \in G$ as an admissible word of length $n$  is called {\it $n$-optimal} if it has the minimum cost among all admissible words of length $n$ that are
equal to $g$.

A decomposition of $g \in G$ as an admissible word of length $n$  is called {\it optimal} if it has the minimum cost among all admissible words of arbitrary lengths that are
equal to $g$.
  
 For a given $g \in G$, an optimal decomposition may not exist since there might be
a sequence of decompositions of $g$ of increasing lengths and with cost going to infinum.
On the other hand, $n$-optimal decompositions exist, as we show in the following
Lemma.

 {\bf Lemma \nopt.} Let $G$ be a connected Lie group, and suppose that the set of controls $S$ is finite and satisfies the following
condition: every generator $C \in S$ with a zero cost factor $\kappa = 0$ has a periodic exponential, i.e., $\exp(T C) = 1$ for some $T > 0$. If $g\in G$ has a decomposition as an admissible word of length $n$ then it has an $n$-optimal decomposition.

{\bf Proof.} Let $\tmin$ be the infinum cost over the set of all decompositions of $g$
in admissible words of length $n$. Let
$$ \left\{ \exp(t_{1}^{(j)} C_{1}^{(j)}) \cdot \ldots \cdot \exp(t_{n}^{(j)} C_{n}^{(j)}) 
\right\}_{j=1,2,\ldots } $$
be a sequence of admissible words of length $n$ with cost converging to $\tmin$.
Since the set of controls $S$ is finite, we can choose a subsequence in which 
the generators $\{ C_{k}^{(j)} \}$ are independent of $j$.

All times  $t_{k}^{(j)}$ are bounded by the same constant (for the generators with non-zero cost factors a bound is obtained from the bound on the total cost, and for other generators from the periodicity assumption). Thus there is a subsequence where all times converge to some values, 
$\lim\limits_{j \to\infty} t_{k}^{(j)} = t_k$.
Then by continuity,
$\exp(t_1 C_1) \cdot \ldots \cdot \exp(t_n C_n) = g$
and the cost of the word $\exp(t_1 C_1) \cdot \ldots \cdot \exp(t_n C_n)$ is $\tmin$.

\

The following Lemma is obvious:

{\bf Lemma \obv.} 
(a) If a word of length $n$ is optimal then it is $n$-optimal.

(b) If  the word 
$$\exp(t_1 C_1) \cdot \ldots \cdot \exp(t_n C_n) \eqno{(\wop)}$$
is optimal (resp. $n$-optimal), then its subword
$$\exp( t_p C_p) \cdot \ldots \cdot \exp(t_k C_k)$$ 
with $1\leq p \leq k \leq n$ is also optimal (resp. $k-p+1$-optimal).   

(c) If the word (\wop) is optimal then
$$\exp(s_1 C_1) \exp(t_2 C_2) \cdot \ldots \cdot \exp(t_{n-1} C_{n-1}) \exp(s_n C_n)$$
with $0\leq s_1\leq t_1$, $0 \leq s_n\leq t_n$, is also optimal.

The last claim of the lemma suggests that for the optimal words there are stronger constraints on time parameters $t_i$ with $2 \leq i \leq n-1$. We will call these the {\it middle time parameters}.

\

{\bf 2. Optimal words in $SU_2$.}

For the rest of the paper we will focus on the case $G = SU_2$. We find it convenient to use the realization of $SU_2$ as the unit sphere in the quaternion algebra $\H$:
$$SU_2 = \left\{ a 1 + bi +cj + dk \, | \, a^2 + b^2 + c^2 + d^2 = 1 \right\} .$$
The Lie algebra $su_2$ in this realization is the tangent space at $1$ and has basis
$\{ i, j, k \}$. This basis is orthonormal relative to the invariant bilinear form.
The norm in $\H$ satisfies $|xy| = |x| \cdot |y|$ and is $SU_2$-bi-invariant. 


The isomorphism with the standard matrix construction of $su_2$ is given
by Pauli matrices:
$$ {i} \mapsto \pmatrix{ i & 0 \cr 0 & -i}, \quad
  {j } \mapsto  \pmatrix{ 0 & 1 \cr -1 & 0}, \quad
  {k} \mapsto  \pmatrix{ 0 & i \cr i & 0} .$$

 The group $SU_2$ acts on its Lie algebra by the conjugation automorphisms. Since $-I$
acts trivially, this action factors through $SO_3 \cong SU_2 / \left\{ \pm I \right\}$,
and is given by the natural action of $SO_3$ on $\R^3$.

The action of $SO_3$ on the unit sphere is transitive, so without loss of generality we
may assume that $X = i$, while $Y = i \cos \alpha + j \sin \alpha$, where $\alpha$ is the angle between the vectors $X$ and $Y$, $0 < \alpha < \pi$. This identification will allow us to carry out certain calculations in an explicit form.
If $C$ is an element of $su_2$ of norm $1$ then $\exp(t C) = \cos(t) + C \sin(t)$.

Since $\exp(\pi X) = \exp(\pi Y) = -1$ is a central element, we see that an
$n$-optimal word satisfies the following 

{\bf $\pi$-Condition:} 
at most one time parameter may be greater or equal to $\pi$; without loss of generality we may assume that this parameter corresponds to the generator $Y$, as it has a lower cost, and is not a middle time parameter.

We begin by describing 4-optimal words.

{\bf Proposition \four.} (a) Let 
$$g = \exp(t_1 X) \exp(t_2 Y) \exp(t_3 X) \exp(t_4 Y) \eqno{(\fourw)}$$
be a 4-optimal word with $t_1, t_2, t_3 < \pi$. Then either
$$ {\tan(t_2) \over \tan(t_3)} = { \kappa - \cos(\alpha) \over 1 - \kappa \cos(\alpha)}
\eqno{(\tanc)}$$
or (\fourw) is not reduced.

(b) The same condition holds for a 4-optimal word
$$\exp(t_4 Y)  \exp(t_3 X)  \exp(t_2 Y) \exp(t_1 X).  $$

{\bf Proof.} We will apply a version of the Lagrange multipliers method. Consider an infinitesimal variation of (\fourw):
$$\exp((t_1 +\ep_1) X) \exp((t_2+\ep_2) Y) \exp((t_3+\ep_3) X) \exp((t_4 + \ep_4)Y) , \eqno{(\fourc)}$$
with the constraint that the product still equals $g$.
Carrying out calculations to the first order in $\ep_i$, we have
$\exp(\ep_i C_i) \approx 1 + \ep_i C_i$. We will collect all terms with $\ep_i$ in the
middle of the word (\fourc), using the relations:
$$\ep_1 X \exp(t_2 Y) = \exp(t_2 Y)  \exp(-t_2 Y)  \ep_1 X \exp(t_2 Y)$$
$$= \exp(t_2 Y)  \left[ (\cos(t_2) - \sin(t_2)Y) \ep_1 X  (\cos(t_2) + \sin(t_2)Y) \right]$$
$$= \exp(t_2 Y)  \ep_1 \left[ \cos^2 (t_2) X - \sin^2(t_2) YXY + \sin (t_2) \cos (t_2)
[X, Y]  \right] .$$
Similarly,
$$\exp(t_3 X) \ep_4 Y = \ep_4 \left[ 
\cos^2 (t_3) Y - \sin^2(t_3) XYX + \sin (t_3) \cos (t_3) [X, Y]  \right] \exp(t_3 X).$$

Setting $Z = {1\over 2} [X, Y] = \sin(\alpha) k$, we can get the following expressions in the basis $\{ X, Y, Z \}:$
$$ XYX =  Y - 2 \cos(\alpha) X, \;\; YXY = X - 2 \cos(\alpha)Y . \eqno{(\XYX)}$$

Using the above relations, we see that to the first order in $\ep_i$, (\fourc) can be written as
$$ \exp(t_1 X) \exp(t_2 Y)
\left[ (\ep_1 \cos(2 t_2) + \ep_3 + 2 \ep_4 \cos(\alpha) \sin^2 (t_3)) X  \right.$$
$$\left. + (2 \ep_1 \cos(\alpha) \sin^2 (t_2) + \ep_2 + \ep_4 \cos(2 t_3)) Y
 + {\kern -1pt} 
(\ep_1 \sin( 2 t_2) 
{\kern -1pt} + {\kern -1pt}
\ep_4 \sin (2 t_3)) Z \right] 
{\kern -3pt}
 \exp(t_3 X) \exp(t_4 Y) .$$

Equating the middle factor to zero, we get the system:
$$\left\{
\matrix{
\ep_1 \cos(2 t_2) + \ep_3 + 2 \ep_4 \cos(\alpha) \sin^2 (t_3) = 0, \cr
2 \ep_1 \cos(\alpha) \sin^2 (t_2) + \ep_2 + \ep_4 \cos(2 t_3) = 0, \cr
\ep_1 \sin( 2 t_2) + \ep_4 \sin (2 t_3)  = 0 . \cr
} \right. 
\eqno{(\fsys)}$$
If this system has rank 3 then by the implicit function theorem there exists a smooth 
curve in the space of parameters $\{ (t_1, t_2, t_3, t_4) \}$ such that the product
(\fourw) is equal identically to $g$ on the curve. The tangent vector to this curve is a non-zero solution of the system (\fsys).  Since by assumption the word (\fourw) is 4-optimal,
the differential of the cost function
$$ \ep_1 + \kappa \ep_2 + \ep_3 + \kappa \ep_4$$
must be zero on the tangent vector of the curve. Thus the determinant
$$\left|
\matrix{
\cos(2 t_2) & 0 & 1 &  2  \cos(\alpha) \sin^2 (t_3) \cr
2 \cos(\alpha) \sin^2 (t_2) & 1 & 0 & \cos(2 t_3)  \cr
\sin( 2 t_2) & 0 & 0 &  \sin (2 t_3)  \cr
1 & \kappa & 1 & \kappa \cr
} \right| $$
must be zero, since the corresponding homogeneous system of equations has a non-trivial
solution. In case when the system (\fsys) has rank less than 3, the above determinant is still equal to zero. Evaluating this determinant we get:
$$ 4 \sin(t_2) \sin(t_3) \left( (1-\kappa \cos(\alpha)) \sin(t_2) \cos(t_3)
- (\kappa - \cos(\alpha)) \cos(t_2) \sin(t_3) \right) = 0.$$
Since $0 < t_2, t_3 < \pi$, we can cancel the factor $ 4 \sin(t_2) \sin(t_3)$ and obtain the claim of the proposition. 

Part (b) is completely analogous.

{\bf Remark \RXX.}  The denominator of  ${ \kappa - \cos(\alpha) \over 1 - \kappa \cos(\alpha)}$ is non-zero since $0 \leq \kappa \leq 1$ and $|\cos(\alpha)| < 1$.
If the numerator of this fraction vanishes, the condition (\tanc) should be replaced with
$$ (1-\kappa \cos(\alpha)) \sin(t_2) \cos(t_3) = (\kappa - \cos(\alpha)) \cos(t_2) \sin(t_3).  \eqno{(\tang)}$$

  Using Proposition \four\ and Lemma \obv\ we get a description of $n$-optimal words:

{\bf Corollary \noptw.} Suppose $ \kappa \neq \cos(\alpha)$. Let 
$ \exp(t_1 C_1) \exp(t_2 C_2)\cdot \ldots \cdot \exp(t_n C_n) $ be a reduced $n$-optimal word with $n \geq 4$.
Then $$t_p = t_{p+2} \hbox{\rm \ for all \ } 2 \leq p \leq n-3. \eqno{(\tae)}$$ 

Since in an $n$-optimal word all middle parameters corresponding to the same control are equal, we will
denote by $\tx$ (resp. $\ty$) the  middle time parameters corresponding to $X$ (resp. $Y$). 

{\bf Corollary \noptx.} Under the assumptions of the previous corollary,
$$  {\tan(\ty) \over \tan(\tx)} =  { \kappa - \cos(\alpha) \over 1 - \kappa \cos(\alpha)}. \eqno{(\tax)}$$
 
 We see from Corollaries \noptw \ and \noptx \ that reduced $n$-optimal words are described with at most three independent time parameters for all $n$. Since the group $SU_2$ is three-dimensional, we conclude that for each $n$ there exists only a finite number of $n$-optimal words representing a given $g \in SU_2$
(for the case $\kappa = \cos(\alpha)$ see Theorem \Tthree\ below).

 Next we shall investigate optimality of words of length $3$. 

{\bf Proposition \triop.} Let $\cos(t) > 0$ and let $\ep>0$ be a small parameter. Then
$$(i) \quad \exp(\ep X) \exp(t Y) \exp(\ep X) = \exp(\tau Y) \exp(\mu X) \exp(\tau Y)
\eqno{(\eXYX)}$$
and
$$(ii) \quad \exp(\ep Y) \exp(t X) \exp(\ep Y) = \exp(\tau X) \exp(\mu Y) \exp(\tau X),
\eqno{(\eYXY)}$$
where
$$\tau = t/2 + \ep \cos(\alpha) (1 - \cos(t)) + o(\ep) \eqno{(\ftau)}$$
and
$$\mu = 2 \ep \cos(t) + o(\ep). \eqno{(\fmu)}$$

{\bf Proof.} Let us write both sides of  (i) in the form $a + b X + c Y + d Z$, using
the relations
$$X Y = - \cos(\alpha) + Z, \quad 
Y X = - \cos(\alpha) - Z, \eqno{(\XY)} $$
together with (\XYX). We immediately see that $d=0$, while
$$a = \cos( 2 \ep) \cos(t) - \cos(\alpha) \sin(2 \ep) \sin(t)
= \cos(2 \tau) \cos(\mu) - \cos(\alpha) \sin(2 \tau) \sin(\mu), \eqno{(\eqa)}$$
$$b  = \sin(2 \ep) \cos(t) - 2\cos(\alpha) \sin^2(\ep) \sin(t)  = \sin(\mu),\eqno{(\eqb)}$$
$$c = \sin(t) = \sin(2 \tau) \cos(\mu) - 2 \cos(\alpha) \sin^2 (\tau)\sin(\mu).\eqno{(\eqc)}$$
We can determine $\mu$ from (\eqb), getting (\fmu).

From the equality 
$$(c+b\cos(\alpha)) \cos(\mu) - a \cos(\alpha) \sin (\mu) = 
\sin(2 \tau)(1 - \sin^2(\alpha) \sin^2(\mu))$$
we determine $\tau$ and obtain (\ftau). It can be seen that this indeed gives a solution of (i). The condition $\cos(t) > 0$ is required to ensure the positivity of $\mu$.
The claim for part (ii) follows by symmetry.

{\bf Corollary \trioptic.}  Let $0 < t < {\pi \over 2}$ and let $\ep>0$ be a small parameter. Then

(i) The word  $\exp(\ep X) \exp(t Y) \exp(\ep X)$ is not 3-optimal.

(ii) If $\kappa > \cos(\alpha)$ then the word
$\exp(\ep Y) \exp(t X) \exp(\ep Y)$ is not 3-optimal.

{\bf Proof.} We are going to show that the cost of the right hand sides of (\eXYX) and (\eYXY) are smaller than for the left hand sides. Let us evaluate the differences in these costs:
$$2 \ep + \kappa t - 2 \kappa \tau - \mu = 
2 \ep + \kappa t -   \kappa t  - 2 \kappa \ep \cos(\alpha) (1 - \cos(t)) - 2 \ep \cos(t) + o(\ep)$$
$$= 2 \ep (1 - \cos(t)) (1 - \kappa \cos(\alpha))  + o(\ep) > 0,$$
which proves (i), while for (ii) we have:
$$t + 2 \kappa \ep - 2 \tau - \kappa \mu
= t + 2 \kappa \ep - t - 2 \ep \cos(\alpha) (1 - \cos(t)) - 2 \kappa \ep \cos(t) + o(\ep)$$
$$= 2 \ep (1 - \cos(t)) (\kappa - \cos(\alpha))  + o(\ep) > 0.$$

\

{\bf Corollary \spi.} Suppose $\kappa > \cos(\alpha)$. A reduced optimal word 
$\exp(t_1 C_1)  \ldots  \exp(t_n C_n)$  of length $n > 1$ 
satisfies the 
{\it strong $\pi$-condition}, which is the $\pi$-condition with an additional restriction $t_1, t_n \leq \pi$.

{\bf Proof.} It is sufficient to show that for a reduced optimal word 
$\exp(t_1 Y) \exp(t_2 X)$ of length 2, we have $t_1 \leq \pi$. Suppose $t_1 > \pi$
and let $0 < \tau < \min(t_2, {\pi\over 2})$.
Then 
$$ \exp(t_1 Y) \exp(t_2 X) = \exp((t_1 - \pi) Y) \exp(\tau X) \exp(\pi Y) \exp((t_2 - \tau) X), $$
which is not optimal by Corollary \trioptic (ii). This contradiction implies $t_1 \leq \pi$.

\

 Now we can describe possible optimal decompositions in $SU_2$. We will consider
several cases.

{\bf Theorem \Tone.} Suppose $\cos(\alpha) < \kappa \leq 1$.  Then the infinum of the cost of admissible decompositions of a given element of $SU_2$ is attained either

(a) on a reduced $n$-optimal word
$$ \quad \exp(t_1 C_1) \exp(t_2 C_2)\cdot \ldots \cdot\exp(t_n C_n) $$
satisfying the strong $\pi$-condition, (\tae) with
${\pi \over 2} \leq \tx < \pi$, ${\pi \over 2} \leq \ty < \pi$, and when $n\geq 4$ the  
condition (\tax). 

\noindent
or

(b)  on a word
$$\exp(t_1 C_1) \exp(t_2 W) \exp(t_3 C_3),$$
where $C_1, C_3 \in \{ X, Y \}$, $t_1, t_2, t_3 \geq 0$  and 
$$W = (1-\kappa \cos(\alpha)) X + (\kappa - \cos(\alpha)) Y$$
with the cost of $\exp(t_2 W)$ equal to $(\kappa^2 - 2 \kappa \cos(\alpha) + 1) t_2$.

{\bf Proof.}  It is clear that the infinum of the cost is either attained on an $n$-optimal word, or is a limit of costs for a sequence of reduced words of length $k$ which are $k$-optimal with $k \to \infty$.

In the first case, the $n$-optimal word must satisfy the conditions in (a) above by
Corollaries \noptw, \noptx, \trioptic\ and \spi.

Alternatively, if for $g \in SU_2$ the infinum of the cost is the limit of costs for a sequence of $k$-optimal words of increasing length, 
$$ \quad g = \exp(t_{1}^{(j)} C_{1}^{(j)}) \exp(t_{2}^{(j)} C_{2}^{(j)})\cdot \ldots \cdot
\exp(t_{{k_j}}^{(j)} C_{{k_j}}^{(j)}), \quad j = 1, 2, 3 \ldots,  \eqno{(\onef)}$$
there will be a subsequence with the fixed first and last generators,
$C_{1}^{(j)} = C^\prime, C_{{k_j}}^{(j)} = C^{\prime\prime}$. 
Since all time parameters belong to the compact set $[0, 2\pi]$, there is a subsequence where $\{ t_{1}^{(j)} \}$
and $\{ t_{{k_j}}^{(j)} \}$ converge,
$$\lim\limits_{j \to \infty} t_{1}^{(j)} = t^\prime, \quad
\lim\limits_{j \to \infty} t_{{k_j}}^{(j)} = t^{\prime\prime}.$$
Time parameters in $k$-optimal words satisfy the conditions of Corollaries \noptw \ and 
\noptx. 
Denote by $\tx^{(j)}$ the middle $X$-times and by $\ty^{(j)}$ the middle $Y$-times in (\onef).
Then $$g = \lim\limits_{j \to \infty}
\exp(t^\prime C^\prime) \left( \exp(\tx^{(j)} X) \exp(\ty^{(j)} Y) \right)^{\left[{k_j \over 2}\right]}
\exp(t^{\prime\prime} C^{\prime\prime}).$$ 
Since the total cost is bounded, we can choose a subsequence for which the sequence
$\left[{k_j \over 2}\right] \tx^{(j)}$ converges, and let 
$$\lim\limits_{j \to \infty} \left[{k_j \over 2}\right] \tx^{(j)} = \tx .$$
By Corollary \noptw
$${\tan(\ty^{(j)}) \over \tan(\tx^{(j)})} = \lambda,
\quad \lambda = { \kappa - \cos(\alpha)  \over  1 - \kappa \cos(\alpha)},$$
and using the fact that 
$\arctan (\lambda \tan(t)) = \lambda t + o(t)$ as $t \to 0$,
we get that  $\lim\limits_{j \to \infty} \left[ {k_j \over 2}\right] \ty^{(j)} = \lambda \tx$ .

By (\prec) we get that 
$$g = \exp(t^\prime C^\prime)  \exp( \tx (X + \lambda Y))  \exp(t^{\prime\prime} C^{\prime\prime}).$$
Finally, we note that $ \tx (X + \lambda Y) = \tilde{t} W$, where 
$W = (1-\kappa \cos(\alpha)) X + (\kappa - \cos(\alpha)) Y$ and 
$\tilde{t} = \tx /  (1-\kappa \cos(\alpha)) $.

{\bf Remark \RXY.} The vector $W$ is orthogonal to the line passing through $X$ and $Y/\kappa$.

{\bf Remark \RXZ.}  Time-optimal decompositions that appear in Theorem \Tone \ involve at most  three independent time parameters. Since the group $SU_2$ is 3-dimensional, there will be a finite number of such decompositions of each length. Moreover, since middle times in the decomposition (a) above are at least ${\pi \over 2}$, any given decomposition of $g$ gives a bound on the length of an optimal decomposition of type (a).

\

 {\bf Theorem \Ttwo.} Suppose $0 < \kappa < \cos(\alpha)$.  Then the infinum of the cost of admissible decompositions of a given element of $SU_2$ is attained on a reduced $n$-optimal word
$$ \quad \exp(t_1 C_1) \exp(t_2 C_2)\cdot \ldots \cdot \exp(t_n C_n) $$
satisfying the  $\pi$-condition, (\tae) with 
 ${\pi \over 2} \leq \ty < \pi$, and when $n\geq 4$ the  
conditions $0 < \tx \leq {\pi \over 2}$ and (\tax).

{\bf Proof.} 
In this case the infinum of the cost can not be realized as a limit cost of a sequence of $n$-optimal words of increasing lengths since by Corollary \noptx\ in an $n$-optimal word we have 
either $\tx \geq {\pi\over 2}$ or $\ty \geq {\pi \over 2}$, 
which implies that the cost would go to infinity as $n \to \infty$.

Applying Proposition \four\ and Corollary \trioptic, we see that for a reduced optimal word of length $n$ the claim of the theorem holds.

\

{\bf Theorem \Tthree.} Suppose $ \kappa = \cos(\alpha) > 0$.  Then the infinum of the cost of admissible decompositions of a given element of $SU_2$ is attained on a reduced $n$-optimal word
$$ \quad \exp(t_1 C_1) \exp(t_2 C_2)\cdot \ldots \cdot \exp(t_n C_n) $$
satisfying the  $\pi$-condition, (\tae) with
 ${\pi \over 2} \leq \ty < \pi$, and when $n\geq 4$ the  
condition $\tx =  {\pi \over 2}$. 

{\bf Proof.} Applying Proposition \four\ (see also Remark \RXX), $n$-optimal words with $n \geq 4$ have middle $X$-times equal to ${\pi \over 2}$. This implies that the infinum of cost is attained on an $n$-optimal word and is not a limit for
a sequence of $n$-optimal words of increasing length. It remains to show that the middle 
$Y$-times are all equal to each other. This can be done using the same method as in the proof of Lemma \trioptic, by considering a the following variation of a 5-optimal word:
$$\exp((t_1 + \ep_1)X) \exp((t_2 + \ep_2) Y)  \exp\left( {\pi \over 2} X\right)
\exp((t_3 + \ep_3)Y) \exp((t_4+ \ep_4)X) .$$
Evaluating the resulting $4 \times 4$ determinant, we get
$$ 8 \cos(\alpha) \sin^2(\alpha) \sin(t_2) \sin(t_3) \sin(t_2 - t_3) = 0,$$
which implies $t_2 = t_3$. 

\

Finally let us consider the case when the cost associated with the generator $Y$ is zero. In this case we have the freedom of replacing the generator $Y$ with $-Y$ since 
$\exp(-tY) = \exp((2\pi -t) Y)$. Thus without loss of generality we may assume that
$\cos(\alpha) \leq 0$. The case when $\cos(\alpha) < 0$ is then covered by Theorem
\Tone.

{\bf Theorem \Tfour.} Let $\cos(\alpha) = 0$ and $\kappa = 0$. For any $g \in SU_2$
the infinum cost is attained on a word of length at most 3.

{\bf Proof.} By Proposition \four,  in an $n$-optimal word of length $n \geq 4$, the middle
$X$-times are equal to ${\pi \over 2}$. Thus the infinum cost is attained on an $n$-optimal word for some $n$. Let for a given $g$, $n$ be the smallest length such that the infinum cost is attained on an $n$-optimal word. 
Suppose the optimal word has a middle factor of $\exp({\pi \over 2} X)$. 
Applying the following identity
$$\exp(t_1 Y) \exp\left({\pi \over 2} X\right) \exp(t_2 Y) = 
\exp((t_1 - t_2) Y) \exp\left({\pi \over 2} X\right) =
\exp\left({\pi \over 2} X\right) \exp((t_2 - t_1)Y) ,$$
we see that such a word is not optimal. This implies $n \leq 3$. 

\

In conclusion, let us consider a particular case.

{\bf Example \Exa.} Let $X = i$, $Y = j$ and $\kappa = 1$. Then an optimal decomposition
of any $g \in SU_2$ is given by words of the following types:

$$(a) \exp( t_1 C_1), \hbox{\rm \ where \ } 0 \leq t_1 < 2\pi; $$

$$(b) \quad \exp(t_1 C_1) \exp( t_2 C_2),
\hbox{\rm \ where \ } 0 < t_1,  t_2 \leq \pi; $$

$$(c) \quad \exp(t_1 C_1) \exp(t_2 C_2) \exp(t_3 C_1), \hbox{\rm \ where \ } 
t_2 \geq {\pi \over 2}, \; t_1, t_2, t_3 \leq \pi; $$ 

$$(d) \quad \exp(t_1 C_1) \exp(t_2 C_2) \exp(t_2 C_1) \exp(t_3 C_2),
 \hbox{\rm \ where \ }  t_2 \geq {\pi \over 2}, \; t_1, t_2, t_3 \leq \pi;$$

or

$$ (e) \quad  \exp(t_1 C_1) \exp\left( t_2 {i+j \over 2}\right) \exp(t_3 C_3),$$
where $C_i \in \{ X, Y \}$, $C_1 \neq C_2$, with the infinum time $\sum_k t_k$.
Here $\exp\left(t {i+j \over 2}\right)$ may be viewed as
$$\exp\left(t {i+j \over 2}\right) = \lim\limits_{N\to\infty} \left[ \exp\left({ti\over 2N}\right)
 \exp\left({tj\over 2N}\right) \right]^N .$$

We obtain this result by applying Theorem \Tone. The only thing we need to show is
that there are no optimal words of length $n \geq 5$. Suppose such a word is indeed optimal. Then it is also $n$-optimal and by Theorem \Tone, all middle times satisfy 
${\pi \over 2} \leq \tx = \ty < \pi$. However,
$$\exp\left({\pi \over 2} i\right) \exp( \ty j) \exp\left({\pi \over 2} i\right) =
 \exp( \pi i)  \exp(- \ty j)= \exp( \pi j) \exp(- \ty j)  = \exp((\pi - \ty) j),$$
which has a lower total time. The analogous equality holds if $i$ and $j$ are switched.
This shows that if the infinum time is attained on a word of length $n$ then $n \leq 4$.
 
 It is fairly straightforward to write down explicit formulas 
for $t_1, t_2, t_3$ in (a)-(e) above for a given $g\in SU_2$.

\

\

\

{\bf References:}

\noindent
[1] N.~Bourbaki,  Lie groups and Lie algebras. Chapters 1-3.  Elements of Mathematics. Springer-Verlag, Berlin, 1989. 

\noindent
[2] D.~Deutsch, A.~Barenco, A.~Ekert,
Universality in quantum computation,
Proc.R.Soc. Lond. A {\bf 449} (1995), 669-677.

\noindent
[3] L.~Euler,
{Formulae generales pro translatione quacunque corporum rigidorum,}
Novi Commentarii Academiae Scientiarum Petropolitanae {\bf 20} (1776), 189-207.

\noindent
[4] V.~Jurdjevic, Geometric control theory, Cambridge studies in advanced mathematics
{\bf 51}, Cambridge University Press, New York, 1997.

\noindent
[5] A.~W.~Knapp, Lie groups beyond an introduction, Progress in Mathematics {\bf 140},
Birkh\"auser, Boston, 1996.

\noindent
[6] S.~Lloyd, Almost any quantum logic gate is universal, 
Phys.Rev.Let. {\bf 75} (1995), 346-349.

\noindent
[7] H.~Rabitz, R.~de Vivie-Riedle, M.~Motzkus, K.~Kompa,
Whither the future of controlling quantum phenomena?
Science {\bf 288} (2000), 824-828.

\end